\documentstyle{article}
\date{}
\title{Hyperbolic Orbits for a Class of Singular Hamiltonian Systems with Repulsive Potentials\footnote{Supported partially by NSF of China.}}
\author{{Donglun Wu and Shiqing Zhang}\\
{\small Yangtze Center of Mathematics and College of Mathematics, Sichuan University,}\\
{\small Chengdu 610064, People's Republic of China}}
\begin{document}
\maketitle \large \baselineskip 14pt

\begin{quote}

{\bf Abstract}\ \ The existence of hyperbolic orbits is proved for a
class of singular Hamiltonian systems with repulsive potentials by
taking limit for a sequence of periodic solutions which are the
minimizers of variational functional.



\end{quote}

\section{Introduction and Main Results}

\ \ \ \ \ \ The motion problem for N bodies with the Newtonian
gravitational force is an old problem since Newton. For the
classical(distinguished from the charged bodies) two-body problem,
Newton and J. Bernoulli have proved that all the moving orbits are
conic curves which depend on the total energy and angular momentum.
As to the classical three-body problems, the mathematicians have not
solved it completely. But if two of the three bodies are restricted
to the moving plane of two-body problem and the third body does not
influence the motion of the first two, this special type of
three-body problem is more simple than the classical one, which is
usually called the restricted three-body problems. In 1987, D. D.
Dionysiou and D. A. Vaiopoulos \cite{24} studied the existence and
locations of the collinear and equilateral Lagrangian points or
solutions for the three-charged-body problems under the rotating
coordinates. When the bodies are charged, they affect each other not
only according to the Newton's, but also to the Coulomb's laws. In
1989, D. D. Dionysiou and G. G. Stamou \cite{25} studied the
stability of motions of the restricted circular and
three-charged-body problems.

In \cite{26}, W. Thirring has interpreted the two-charged-body
problems. The motion equations are the following second order
Hamiltonian systems
\begin{equation}
   \ddot{u}(t)+\nabla V(u(t))=0,\label{1}
\end{equation}
where the kinetic energy $K$ and potential $V$ have the form of
\begin{eqnarray*}
K=\frac{|p_{1}|^{2}}{2m_{1}}+\frac{|p_{2}|^{2}}{2m_{2}},\ \ \ \ \ \
V=\frac{\delta}{|x_{1}-x_{2}|}
\end{eqnarray*}
with $\delta=e_{1}e_{2}-\kappa m_{1}m_{2}$, where $p_{i}$ represent
the momentum; $e_{i}$ represent the charge; $x_{i}$ represent the
position and $m_{i}$ represent the mass of the $i$th$(i=1,2)$ body;
$\kappa$ is the gravitational constant. The energy function is
\begin{eqnarray*}
H=\frac{|p_{1}|^{2}}{2m_{1}}+\frac{|p_{2}|^{2}}{2m_{2}}+\frac{\beta}{|x_{1}-x_{2}|}.
\end{eqnarray*}

It is well known that, by some simple calculations, we can separate
two-body problems to the center-of-mass and the relative
coordinates. Furthermore, if we choose the relative coordinates, we
can reduce the kinetic energy and the potential energy to the
following forms
\begin{eqnarray*}
\overline{K}=\frac{|p|^{2}(m_{1}+m_{2})}{2m_{1}m_{2}},\ \ \ \ \ \
\overline{V}=\frac{\beta}{|x|},
\end{eqnarray*}
where
\begin{eqnarray*}
p=\frac{m_{1}p_{1}}{m_{1}+m_{2}}-\frac{m_{2}p_{2}}{m_{1}+m_{2}}\ \ \
\mbox{and}\ \ \ \ x=x_{1}-x_{2}.
\end{eqnarray*}
Moreover, the energy function takes the following form of
\begin{eqnarray*}
H=\frac{|p_{cm}|^{2}}{2(m_{1}+m_{2})}+\frac{|p|^{2}(m_{1}+m_{2})}{2m_{1}m_{2}}+\frac{\beta}{|x|}\doteq
H_{s}+H_{r},
\end{eqnarray*}
where $p_{cm}=p_{1}+p_{2}$ and $H_{r}$ is a limiting case of $H$, in
which one body has infinite mass and the other has the reduce mass,
which is
\begin{eqnarray*}
m=\frac{(m_{1}+m_{2})}{2m_{1}m_{2}}.
\end{eqnarray*}
With the developments of variational methods, more and more
mathematicians use the variational methods to look for the periodic,
homoclinic and heteroclinic orbits of Hamiltonian systems
[1-3,5,10-18,20-23,25,27-30] and the references therein. But for the
classical two-body problems, there are only a few papers involving
the existence of hyperbolic orbits via the variational methods with
fixed energy. A main difficulty is how to prove the obtained orbit
is not at the infinite point. If the charges of the bodies take the
same sign(both positive or both negative) and the charges are large
enough, we can see that $\delta>0$, which is much different from the
classical Newtonian two-body problems, since the potential is
positive and the effect between the bodies becomes repulsion rather
than attraction.

In Thirring's statement, the projection of trajectory onto
$R^{3}\times R^{3}\setminus\{x:x=0\}$ under the polar coordinates is
\begin{eqnarray*}
r=\frac{|L|^{2}}{|F|\cos\zeta-m\delta},
\end{eqnarray*}
where $\zeta=\angle(F,x)$ and $L$, $F$ are constants of the motion
stated as follows
\begin{eqnarray*}
L=[x\times p],\ \ \ \ \ \ \ F=[p\times L]+m\delta\frac{x}{|x|}.
\end{eqnarray*}
Actually, $L$ is the angular momentum and $F$ is known as Lenz
vector. The relationship among $L$, $F$ and $H_{r}$ is
\begin{eqnarray*}
|F|^{2}=2m|L|^{2}H_{r}+m^{2}\delta^{2}.
\end{eqnarray*}
Obviously, $H_{r}>0$ leads to that $|F|>|m\delta|$ and $r$ becomes
infinite at $\zeta=\arccos(m\delta/|F|)$. Then the trajectory is
hyperbolic(or linear, if $\delta=0$). In this paper, an orbit of
(\ref{1}) is said to be hyperbolic, if it satisfies
\begin{eqnarray*}
|u(t)|\rightarrow\infty\ \ \ \mbox{as}\ \ \
t\rightarrow\pm\infty.\label{5}
\end{eqnarray*}

In the present paper, we are concerned with the repulsive potential
which is positive(when potential is negative, it is very similar to
the classical case) and $(-\alpha)$-homogene-ous($0<\alpha<2$, which
equals to 1 in the classical models). More precisely, we consider
the system (\ref{1}) with
\begin{eqnarray}
   \frac{1}{2}|\dot{u}(t)|^{2}+ V(u(t))=H\label{2},
\end{eqnarray}
where $u\in(R^{1},R^{N})$, $V(>0)\in
C^{1}(R^{N}\setminus\{0\},R^{1})$ has a singularity at the origin.
Subsequently, $\nabla V(x)$ denotes the gradient with respect to the
$x$ variable, $(\cdot,\cdot):R^{N}\times R^{N} \rightarrow  R$\
denotes the standard Euclidean inner product in $R^{N}$ and $\mid
\cdot\mid$ is the induced norm. In 2000, for $N=2$, Felmer and
Tanaka \cite{3} proved that

\vspace{0.3cm} \textbf{Theorem 1.1(See\cite{3}).}\ {\em Assume that $N=2$ and the following conditions hold\\

$(A_{1})$ $V\in C^{1}(R^{N}\setminus\{0\},R^{1})$,

$(A_{2})$ $V(x)\leq0$ for all $x\in R^{N}\setminus\{0\}$,

$(A_{3})$ there are constants $\tau>2$, $r_{0}>0$ and $d_{0}>0$ such
that
\begin{eqnarray*}
&&(i).\ -V(x)\geq \frac{d_{0}}{|x|^{\tau}}\ \ \ \mbox{for}\ \ 0<|x|\leq r_{0},\nonumber\\
&&(ii).\ (x,\nabla V(x))+2V(x)\rightarrow +\infty\ \mbox{as}\
|x|\rightarrow0.\nonumber
\end{eqnarray*}

$(A_{4})$  there exist $\lambda>1$ and $k_{0}>0$ such that
\begin{eqnarray*}
-V(x)\leq\frac{k_{0}}{|x|^{\lambda}}\ \ \ \mbox{and}\ \ \ |\nabla
V(x)|\leq\frac{k_{0}}{|x|^{\lambda+1}}\ \ \ \mbox{for}\ \ |x|\geq1.
\end{eqnarray*}

Then for any given $H>0$, $\gamma_{+}$, $\gamma_{-}\in R$ with
$\gamma_{+}-\gamma_{-}>\pi$, there exists a solution
$u(t)=r(t)(\cos\gamma(t),\sin\gamma(t))$ of $(1)-(2)$such that
$\gamma\rightarrow\gamma_{\pm}$ as
$t\rightarrow\pm\infty$.}\\

For $N\geq3$, they proved that

\vspace{0.3cm}{\bf Theorem 1.2(See\cite{3})}\ \ {\em Assume $N\geq3$
and $(A_{1})-(A_{4})$ hold. Then for any given $H>0$ and
$\gamma_{+}\neq-\gamma_{-}$, there exists a solution $u(t)$ of
$(1)-(2)$ such that}
\begin{eqnarray*}
\lim_{t\rightarrow\pm\infty}\frac{u(t)}{|u(t)|}=\gamma_{\pm},
\end{eqnarray*}
where $\gamma_{+}$, $\gamma_{-}\in S^{N-1}=\{x\in R^{N}||x|=1\}$ are
the asymptotic directions for the solution $u(t)$.

In 2011, Wu and Zhang in \cite{5} proved the existence of the
hyperbolic orbits for another class of singular Hamiltonian systems.
They obtained the following theorem.

\vspace{0.3cm}{\bf Theorem 1.3(See\cite{5})}\ \ {\em Suppose that
$V\in C^{1}(R^{N}\setminus\{0\},R^{1})$ satisfies

\vspace{0.3cm}{\bf $(B_{1})$}\ $V(-x)=V(x)$, $\forall x\in
R^{N}\setminus\{0\}$,

\vspace{0.3cm}{\bf $(B_{2})$}\ $V(x)<0$, $\forall x\in
R^{N}\setminus\{0\}$,

\vspace{0.3cm}{\bf $(B_{3})$}\ $2 V(x)+(x,\nabla V(x))\rightarrow0\
\ \ \mbox{as}\ \ \ |x|\rightarrow+\infty$,

\vspace{0.3cm}{\bf $(B_{4})$}\ $2 V(x)+(x,\nabla
V(x))\rightarrow+\infty\ \ \ \mbox{as}\ \ \ |x|\rightarrow0$,

\vspace{0.3cm}{\bf $(B_{5})$}\ $-V(x)\rightarrow+\infty$ as
$|x|\rightarrow0$,

\vspace{0.3cm}{\bf $(B_{6})$}\ $V(x)\rightarrow0$ as
$|x|\rightarrow+\infty$.
\\

Then for any $H>0$, there is at least one hyperbolic orbit for
systems  $(1)-(2)$.}

Motivated by the above papers, we study systems $(1)-(2)$, under
some assumptions, we obtain the hyperbolic orbits
 with $H>0$. Precisely, we prove the
following theorem.

\vspace{0.3cm}{\bf Theorem 1.4}\ \ {\em Suppose that $V\in
C^{1}(R^{N}\setminus\{0\},R^{1})$ satisfies $(B_{1})$ and

\vspace{0.3cm}{\bf $(V_{1})$}\ there is a constant $\alpha\in (0,2)$
such that
\begin{eqnarray*}
(x,\nabla V(x))= -\alpha V(x)<0\ \ \ \mbox{for}\ \ \mbox{any}\ \ \
x\in R^{N}\setminus\{0\}.
\end{eqnarray*}

Then for any $H>0$, there is at least one hyperbolic orbit for
systems  $(1)-(2)$.}

\vspace{0.3cm}{\bf Remark 1}\ It is easy to see that
$V(x)=\displaystyle\frac{1}{|x|^{\alpha}}(0<\alpha<2)$ satisfies our
conditions. It is well known that, for the potential
$V(x)=-\displaystyle\frac{1}{|x|^{\alpha}}$, $\alpha\geq2$ is called
the strong force case; $0<\alpha<2$ is called the weak force case.
The strong force case was firstly studied by
Poincar$\acute{\mbox{e}}$ in 1896 to avoid the minimizing sequence
converging to some collision point. Obviously, Theorem 1.1-Theorem
1.3 treat the strong force case. $(B_{4})$ is another form for the
classical strong force conditions. As to the weak force cases,
mathematicians usually perturbed the potential such that it
satisfies the so-called $Gordon's\ Strong\ Force$ condition which is
introduced by Gordon\cite{7} in 1975.

Under some additional conditions, we can get the asymptotic
direction of the solution at infinity. We have the following
theorem.

\vspace{0.3cm}{\bf Theorem 1.5}\ \ {\em Suppose that $V\in
C^{1}(R^{N}\setminus\{0\},R^{1})$ satisfies $(B_{1})$, $(V_{1})$ and
the following condition

\vspace{0.3cm}{\bf $(V_{2})$}\ there exist constants $\beta>1$,
$M_{0}>0$ and $r_{0}\geq1$  such that
\begin{eqnarray*}
|x|^{\beta+1}|\nabla V(x)|\leq M_{0},\ \ \ \mbox{for all}\ \ \
|x|\geq r_{0}.
\end{eqnarray*}

Then for any $H>0$, there is at least one hyperbolic orbit for
systems  $(1)-(2)$ which has the given asymptotic direction at
infinity.}

\vspace{0.3cm}{\bf Remark 2}\ It is easy to see that
$V(x)=\displaystyle\frac{1}{|x|^{\beta}}(1<\beta<2)$ satisfies our
conditions.

In this paper, the potential is a sum gotten simultaneously
according to Newton's and Coulomb's laws, but the sum is positive
which is different from the classical gravitational case, the
negative potential is out of our study. But the following remark
shows the connection.

\vspace{0.3cm}{\bf Remark 3}\ \  Suppose that $V$ satisfies the
 condition $(V_{1})$, then $V$ satisfies $(B_{3})$ and
 $(B_{4})$. The proof can easily be obtained from Lemma 2.1 in the following section.

\section{Variational Settings}

\ \ \ \ \ \ Let us set
\begin{eqnarray*}
&&H^{1}=W^{1,2}(R^{1}/Z,R^{N}),\nonumber\\
&&E=\{q\in H^{1}|\ q(t)\neq0,\forall t\in [0,1]\}.
\end{eqnarray*}

Moreover, let $L^{\infty}([0,1],R^{N})$ be a space of measurable
functions from $[0,1]$ into $R^{N}$ and essentially bounded under
the following norm
\begin{eqnarray*}
\ \|q\|_{L^{\infty}([0,1],R^{N})}:=esssup\{|q(t)|:t\in[0,1]\}.
\end{eqnarray*}
Let $f$: $E\rightarrow R^{1}$ be the functional defined by
\begin{eqnarray*}
f(q)&=&\frac{1}{2}\int^{1}_{0}|\dot{q}(t)|^2dt\int^{1}_{0}(H-V(q(t)))dt.
\end{eqnarray*}
Then one can easily check that
\begin{eqnarray*}
\ \langle
f'(q),q\rangle&=&\int^{1}_{0}|\dot{q}(t)|^2dt\int^{1}_{0}\left(H-V(q(t))-\frac{1}{2}(
\nabla V(q(t)),q(t))\right)dt.
\end{eqnarray*}

As Tanaka stated in \cite{27}, when the potential
$V=-\frac{1}{|x|^{\alpha}}(\alpha>0)$, we have

\vspace{0.3cm}{\bf Case 1}\ $\alpha>2$(strong force).

(i) For $H>0$, system (\ref{1}) - (\ref{2}) possesses an explicit
solution
\begin{eqnarray*}
q(t)=r_{\alpha}(\cos\omega_{\alpha}t,\sin\omega_{\alpha}t,0,\cdots,0),
\end{eqnarray*}
where $r_{\alpha}=(\frac{\alpha-2}{2H})^{\frac{1}{\alpha}}$ and
$\omega_{\alpha}=(\alpha r_{\alpha}^{-(\alpha+2)})^{\frac{1}{2}}$.

(ii) For $H\leq0$, we have $\langle f'(q),q\rangle>0$ for all
non-constant $q\in E$. Thus system (\ref{1}) - (\ref{2}) possesses
no periodic solutions.

\vspace{0.3cm}{\bf Case 2}\ $\alpha\in(0,2)$(weak force).

(i) For $H\geq0$, we have $\langle f'(q),q\rangle>0$ for all
non-constant $q\in E$. Thus system (\ref{1}) - (\ref{2}) possesses
no periodic solutions.

(ii) For $H<0$, system (\ref{1}) - (\ref{2}) possesses an explicit
solution
\begin{eqnarray*}
q(t)=r_{\alpha}(\cos\omega_{\alpha}t,\sin\omega_{\alpha}t,0,\cdots,0),
\end{eqnarray*}
where $r_{\alpha}=(\frac{\alpha-2}{2H})^{\frac{1}{\alpha}}$ and
$\omega_{\alpha}=(\alpha r_{\alpha}^{-(\alpha+2)})^{\frac{1}{2}}$.

But in our model, the potential and the total energy are both
positive. Similar to A. Ambrosetti and V. Coti Zelati in \cite{1}
and some early papers \cite{28,30,29,31}, we consider the following
set
\begin{eqnarray*}
M_{H}=\left\{q\in E| \int^{1}_{0}(V(q(t))+\frac{1}{2}(\nabla
V(q(t)),q(t)))dt=H\right\}
\end{eqnarray*}
and for any given unite vector(direction) $e\in S^{N-1}$, we set
\begin{eqnarray*}
\Lambda_{R}=\{q\in M_{H}|\ q(t+\frac{1}{2})=-q(t),q(0)=q(1)=Re\}.
\end{eqnarray*}
We fix the direction of the vectors at time 0 and 1 in $\Lambda_{R}$
to prove the asymptotic directions at infinity of the hyperbolic
orbits we obtained can be the same one. For any $q\in H^{1}$, we
know that the following norms are equivalent to each other
\begin{eqnarray}
&&\|q\|_{H^1}=\left(\int^{1}_{0}|\dot{q}(t)|^2dt\right)^{1/2}+\left|\int^{1}_{0}q(t)dt\right|\nonumber\\
&&\|q\|_{H^1}=\left(\int^{1}_{0}|\dot{q}(t)|^2dt\right)^{1/2}+\left(\int^{1}_{0}|q(t)|^2dt\right)^{1/2}\nonumber\\
&&\|q\|_{H^1}=\left(\int^{1}_{0}|\dot{q}(t)|^2dt\right)^{1/2}+|q(0)|.\nonumber
\end{eqnarray}

If $q\in \Lambda_{R}$, we have $\displaystyle\int^{1}_{0}q(t)dt=0$,
then by Poincar$\acute{\mbox{e}}$-Wirtinger's inequality, we obtain
that the above norms are equivalent to
\begin{eqnarray*}
\|q\|_{H^1}=\left(\int^{1}_{0}|\dot{q}(t)|^2dt\right)^{1/2}.
\end{eqnarray*}

We remark that under the following condition
\begin{eqnarray}
3(x,\nabla V(x))+(x,\nabla^{2} V(x)x)\neq0, \ \ \ \ \mbox{for all} \
\ \ x\in R^{N}\setminus\{0\},\label{4}
\end{eqnarray}
A. Ambrosetti and V. Coti Zelati have proved that $M_{H}$ is a
non-empty $C^{1}$ manifold of codimension 1 in $E$ and all critical
points of $f$ on $E$ belongs to $M_{H}$. They develop the
Ljusternik-Schnirelman theory for the restricted functional $f$:
$M_{H}\rightarrow R^{1}$ to find critical points. By the setting of
$M_{H}$, we study the following functional $f$
\begin{eqnarray*}
f(q)&=&\frac{1}{4}\|q\|^{2}\int^{1}_{0}(\nabla V(q(t)),q(t)))dt,\ \
\ \ q\in \Lambda_{R}.
\end{eqnarray*}

In this paper, we need to fix the end point such that
$|q(0)|=|q(1)|=R$ and restrict the symmetry of the functions which
is different from the case in \cite{1}. Firstly, we prove the set
$\Lambda_{R}$ is not empty. In order to do this, we need the
following lemma.

\vspace{0.3cm}{\bf Lemma 2.1}\ \ {\em Suppose that $(V_{1})$ holds.
Then we have}
\begin{eqnarray*}
\ C_{1}|x|^{-\alpha}\leq V(x)\leq C_{2}|x|^{-\alpha}\ \ \
\mbox{for}\ \ x \in R^{N}\setminus\{0\},\label{24}
\end{eqnarray*}
where $C_{1} \stackrel{\Delta}{=}\inf\{V(x)|\ x \in R^{N},|x|=1\}$,
$C_{2} \stackrel{\Delta}{=}\sup\{V(x)|\ x \in R^{N} ,|x|=1\}$.

\vspace{0.3cm}{\bf Proof.}\ \ Set $\varphi(s)=s^{\alpha}V(s\theta)$.
By
 $(V_{1})$, we have
\begin{eqnarray*}
\ \varphi^{'}(s)&=&s^{\alpha}(\nabla V(s\theta),\theta)+\alpha s^{\alpha-1}V(s\theta)\\
&=& \mbox{} s^{\alpha-1}\left((\nabla V(s\theta),s\theta)+\alpha V(s\theta)\right)\\
&=& \mbox{} 0,
\end{eqnarray*}
then, for any $s\in R^{1}\setminus\{0\}$, we have
\begin{eqnarray*}
\ s^{\alpha}V(s\theta)= V(\theta).
\end{eqnarray*}
Set $s=|x|$ and $\theta=x/|x|$, we obtain that
\begin{eqnarray*}
V(x)= V\left(\frac{x}{|x|}\right)|x|^{-\alpha}\ \ \ \mbox{for}\ \ x
\in R^{N}\setminus\{0\},
\end{eqnarray*}
which implies that
\begin{eqnarray*}
V(x)\leq C_{2}|x|^{-\alpha}\ \ \ \mbox{for}\ \ x \in
R^{N}\setminus\{0\}
\end{eqnarray*}
and
\begin{eqnarray*}
\ C_{1}|x|^{-\alpha}\leq V(x)\ \ \ \mbox{for}\ \ x \in
R^{N}\setminus\{0\},
\end{eqnarray*}
which proves the lemma.

\vspace{0.3cm}{\bf Remark 4}\ From $(V_{1})$ and Lemma 2.1, we can
easily prove that $V$ satisfies hypotheses $(B_{3})$, $(B_{4})$ and
$(B_{6})$ in Theorem 1.3.

The following lemma shows that $\Lambda_{R}$ is not empty.

\vspace{0.3cm}{\bf Lemma 2.2}\ \ {\em Suppose that $(V_{1})$ holds
and let
\begin{eqnarray*}
g(q)=\displaystyle\int^{1}_{0}(V(q(t))+\frac{1}{2}(\nabla
V(q(t)),q(t)))dt.
\end{eqnarray*}
Then, for any $H>0$, the equation $g(q)=H$ has at least one solution
on $E$ such that $q(t+\frac{1}{2})=-q(t), q(0)=q(1)=Re$.}

\vspace{0.3cm}{\bf Proof.}\ \ For any fixed total energy $H$, we set
$q_{a}(t)=Re\cos^{p}(2\pi t)+a\eta\sin^{p}(2\pi t)$, $t\in [0,1]$,
where $\eta \in R^{N}$, $|\eta|=1$, $\langle e,\eta\rangle=0$,
$a>0$, $p=2\left[\frac{2}{\alpha}\right]+1$. It is obvious that
$q_{a}\in E$ and $q_{a}(t+\frac{1}{2})=-q_{a}(t),q(0)=q(1)=Re$. It
follows from  $(V_{1})$ that
\begin{eqnarray*}
g(q_{a})=(1-\displaystyle\frac{\alpha}{2})\displaystyle\int^{1}_{0}V(q_{a}(t))dt>0.
\end{eqnarray*}
By Lemma 2.1, we can deduce that
\begin{eqnarray}
\frac{C_{1}(2-\alpha)}{2}\int^{1}_{0}\frac{1}{|q_{a}(t)|^{\alpha}}dt\leq
g(q_{a})\leq
\frac{C_{2}(2-\alpha)}{2}\int^{1}_{0}\frac{1}{|q_{a}(t)|^{\alpha}}dt.\label{6}
\end{eqnarray}

Let $I_{1}=\{t\in[0,1]| \sin^{2p}(2\pi t)\leq a^{-1}\}$,
$I_{2}=[0,1]\setminus I_{1}$, which implies that
$|I_{1}|=\displaystyle\frac{2\arcsin a^{-\frac{1}{2p}}}{\pi}$. Then
we have

\begin{eqnarray*}
& & \int^{1}_{0}\frac{1}{|q_{a}(t)|^{\alpha}}dt\nonumber\\&=&
\mbox{}\int^{1}_{0}\frac{1}{(R^{2}\cos^{2p}(2\pi
t)+a^{2}\sin^{2p}(2\pi t))^{\frac{\alpha}{2}}}dt\nonumber\\
&=& \mbox{} \int_{I_{1}}\frac{1}{(R^{2}\cos^{2p}(2\pi
t)+a^{2}\sin^{2p}(2\pi
t))^{\frac{\alpha}{2}}}dt+\int_{I_{2}}\frac{1}{(R^{2}\cos^{2p}(2\pi
t)+a^{2}\sin^{2p}(2\pi
t))^{\frac{\alpha}{2}}}dt\nonumber\\
&\leq& \mbox{} \int_{I_{1}}\frac{1}{(R^{2}\cos^{2p}(2\pi
t))^{\frac{\alpha}{2}}}dt+\int_{I_{2}}\frac{1}{(a^{2}\sin^{2p}(2\pi
t))^{\frac{\alpha}{2}}}dt\nonumber\\
&\leq& \mbox{}
\frac{|I_{1}|}{R^{\alpha}(1-a^{-\frac{1}{p}})^{\frac{p\alpha}{2}}}+a^{-\frac{\alpha}{2}}\nonumber\\
&\leq& \mbox{} \frac{2\arcsin a^{-\frac{1}{2p}}}{\pi
R^{\alpha}(1-a^{-\frac{1}{p}})^{\frac{p\alpha}{2}}}+a^{-\frac{\alpha}{2}},
\end{eqnarray*}
which implies that for any $R>0$, we have
\begin{eqnarray*}
\int^{1}_{0}\frac{1}{|q_{a}(t)|^{\alpha}}dt\rightarrow0\ \ \
\mbox{as}\ \ \ a\rightarrow+\infty.
\end{eqnarray*}
It follows from (\ref{6}) that
\begin{eqnarray}
g(q_{a})\rightarrow0\ \ \ \mbox{as}\ \ \
a\rightarrow+\infty.\label{8}
\end{eqnarray}
On the other hand, let $I_{3}=\{t\in[0,1]| \cos^{2p}(2\pi t)\leq
a\}$. Then we obtain
\begin{eqnarray*}
|I_{3}|=\left(1-\frac{2}{\pi}\arccos\left(a^{\frac{1}{2p}}\right)\right)\
\ \ \ \mbox{and}\ \ \ \
\frac{|I_{3}|}{a^{\frac{1}{2p}}}\rightarrow\frac{2}{\pi}\ \ \
\mbox{as}\ \ \ a\rightarrow0.
\end{eqnarray*}
For any $a>0$, it follows from the definition of $q_{a}$ that
\begin{eqnarray*}
\int^{1}_{0}\frac{1}{|q_{a}(t)|^{\alpha}}dt&\geq&
\int_{I_{3}}\frac{1}{(R^{2}\cos^{2p}(2\pi t)+a^{2}\sin^{2p}(2\pi
t))^{\frac{\alpha}{2}}}dt\nonumber\\
&\geq& \mbox{}
\frac{|I_{3}|}{(R^{2}a+a^{2})^{\frac{\alpha}{2}}}\nonumber\\
&\geq& \mbox{}
\frac{|I_{3}|}{a^{\frac{1}{2p}}}\frac{a^{\frac{1}{2p}}}{(R^{2}a+a^{2})^{\frac{\alpha}{2}}}\nonumber\\
&\geq& \mbox{}
\frac{|I_{3}|}{a^{\frac{1}{2p}}}\frac{1}{\left(R^{2}a^{\left(1-\frac{1}{\alpha
p}\right)}+a^{\left(2-\frac{1}{\alpha
p}\right)}\right)^{\frac{\alpha}{2}}}.\label{48}
\end{eqnarray*}
By the definition of $p$ and $\alpha\in(0,2)$, we can see that
$2-\displaystyle\frac{1}{\alpha p}>1-\displaystyle\frac{1}{\alpha
p}>0$. Then we can deduce that
\begin{eqnarray*}
\left(R^{2}a^{\left(1-\frac{1}{\alpha
p}\right)}+a^{\left(2-\frac{1}{\alpha
p}\right)}\right)^{\frac{\alpha}{2}}\rightarrow0\ \ \ \mbox{as}\ \
a\rightarrow0,
\end{eqnarray*}
with $\displaystyle\frac{|I_{3}|}{a^{\frac{1}{2p}}}>\frac{1}{\pi}$
when $a$ is near $0$, which implies that
\begin{eqnarray*}
\int^{1}_{0}\frac{1}{|q_{a}(t)|^{\alpha}}dt\rightarrow+\infty\ \ \
\mbox{as}\ \ a\rightarrow0.
\end{eqnarray*}
It follows from (\ref{6}) that
\begin{eqnarray}
g(q_{a})\rightarrow+\infty\ \ \ \mbox{as}\ \ \
a\rightarrow0.\label{9}
\end{eqnarray}
Combining (\ref{8}) and (\ref{9}), we obtain that equation $g(q)=H$
has at least one solution in $E$, for any $H>0$, such that
$q(t+\frac{1}{2})=-q(t),q(0)=q(1)=Re$, which implies that
$\Lambda_{R}\neq \emptyset.$ The proof of this lemma is completed.

For any $q\in M_{H}$, by the homogeneous property of $V$, we have
\begin{eqnarray*}
\int^{1}_{0}(H-V(q(t)))dt=\frac{1}{2}\int^{1}_{0}(\nabla
V(q(t)),q(t)))dt=-\frac{\alpha}{2}\int^{1}_{0}V(q(t))dt,
\end{eqnarray*}
which implies
\begin{eqnarray}
\int^{1}_{0}V(q(t))dt=\frac{2H}{2-\alpha}.\label{10}
\end{eqnarray}
In our model, since the potential and the total energy are both
positive, if we want to use the minimizing theory to get the
critical points which yield the periodic solutions of system
(\ref{1}) - (\ref{2}), we need to modify our functional as follow.
\begin{eqnarray}
F(q)=-f(q)&=&-\frac{1}{4}\|q\|^{2}\int^{1}_{0}(\nabla
V(q(t)),q(t)))dt\nonumber\\
&=& \mbox{}\frac{\alpha}{4}\|q\|^{2}\int^{1}_{0}V(q(t))dt>0.
\end{eqnarray}
It is easy to see that $F$ and $f$ share the same critical points.
Our way to get the hyperbolic orbit is by approaching it with a
sequence of periodic solutions. The approximate solutions are
obtained by the minimizing theory. We need the following lemma which
is proved by A. Ambrosetti and V. Coti Zelati in \cite{1}.

\vspace{0.3cm}{\bf Lemma 2.3(See\cite{1})}\ {\em Let
$f(q)=\displaystyle\frac{1}{2}\int^{1}_{0}|\dot{q}(t)|^2dt\int^{1}_{0}(H-V(q(t)))dt$
and $\tilde{q}\in H^{1}$ be such that $f^{'}(\tilde{q})=0$,
$f(\tilde{q})>0$. Set
\begin{eqnarray*}
T^{2}=\frac{\displaystyle\frac{1}{2}\displaystyle\int^{1}_{0}|\dot{\tilde{q}}(t)|^{2}dt}{\displaystyle\int^{1}_{0}(H-V(\tilde{q}(t))dt}.
\end{eqnarray*}
Then $\tilde{u}(t)=\tilde{q}(t/T)$ is a non-constant $T$-periodic
solution for (\ref{1}) and (\ref{2}).}

\vspace{0.3cm}{\bf Remark 5}\ In view of the proof of Lemma 2.3 in
\cite{1}, we can see that the condition $f(\tilde{q})>0$ in Lemma
2.3 can be replaced by
$\displaystyle\int^{1}_{0}|\dot{\tilde{q}}(t)|^{2}dt>0$.

\vspace{0.3cm}{\bf Lemma 2.4(Palais\cite{13})}\ {\em Let $\sigma$ be
an orthogonal representation of a finite or compact group $G$ in the
real Hilbert space $H$ such that for any $\sigma\in G$,
\begin{eqnarray*}
f(\sigma\cdot x)=f(x),
\end{eqnarray*}
where $f\in C^{1}(H,R^{1})$. Let $S=\{x\in H|\sigma x=x,
\forall\sigma\in G\}$, then the critical point of $f$ in $S$ is also
a critical point of $f$ in $H$.}

\vspace{0.3cm}{\bf Lemma 2.5(Translation Property\cite{8})}\ {\em
Suppose that, in domain $D\subset R^{N}$, we have a solution
$\phi(t)$ for the following differential equation
\begin{eqnarray*}
x^{(n)}+F(x^{(n-1)},\cdots,x)=0,
\end{eqnarray*}
where $x^{(k)}=d^{k}x/dt^{k}$, $k=0,1,\cdots,n$, $x^{(0)}=x$. Then
$\phi(t-t_{0})$ with $t_{0}$ being a constant is also a solution.}

Firstly, we prove the existence of the approximate solutions, then
we study the limit procedure.

\section{Existence of Periodic Solutions}

\ \ \ \ \ \ \ In order to obtain the critical points of the
functional and make some estimations, we need the following lemma.

\vspace{0.3cm}{\bf Lemma 3.1}\ {\em Suppose the conditions of
Theorem 1.4 hold, then for any  $R>0$, there exists at least one
periodic solution on $\overline{\Lambda}_{R}$ for the following
systems
\begin{equation}
   \ddot{q}(t)+\nabla V(q(t))=0,\ \ \ \
   \forall\ t\in\left(-\frac{T_{R}}{2},\frac{T_{R}}{2}\right)\label{16}
\end{equation}
with
\begin{equation}
   \frac{1}{2}|\dot{q}(t)|^{2}+ V(q(t))=H,\ \ \ \ \ \
   \forall\ t\in\left(-\frac{T_{R}}{2},\frac{T_{R}}{2}\right).\label{17}
\end{equation}}

\vspace{0.3cm}{\bf Proof.}\ We notice that $H^{1}$ is a reflexive
Banach space and $\overline{\Lambda}_{R}$ is a weakly closed subset
of $H^{1}$. By the definition of $F$ and (\ref{10}), we obtain that
$F$ is a functional bounded from below and
\begin{eqnarray}
F(q)&=&\frac{\alpha}{4}\|q\|^{2}\int^{1}_{0}V(q(t))dt\nonumber\\
&=& \mbox{}\frac{\alpha H}{2(2-\alpha)}\|q\|^{2}\rightarrow +\infty\
\ \ \mbox{as}\ \ \|q\|\rightarrow +\infty.
\end{eqnarray}
Furthermore, it is easy to check that $F$ is weakly lower
semi-continuous. Then, we can see that for every $R>0$ there exists
a minimizer $q_{R}\in \overline{\Lambda}_{R}$ such that
\begin{eqnarray}
F'(q_R)=0,\ \ \ \ F(q_R)=\inf_{q\in \overline{\Lambda}_{R}}
F(q)\geq0.\label{18}
\end{eqnarray}
It is easy to see that
$\|q\|^{2}=\int^{1}_{0}|\dot{q}_{R}(t)|^{2}dt>0$, otherwise we
deduce that $|q_{R}(t)|\equiv R$, on the other hand, by the
$1/2$-antisymmetry of $q_{R}$, we have $q_{R}\equiv 0$, which is a
contradiction. This implies that $F(q_R)>0$. By the definition of
$F$, we have
\begin{eqnarray}
f'(q_R)=0.
\end{eqnarray}

Then let
\begin{eqnarray}
T^{2}_{R}=\frac{\displaystyle\frac{1}{2}\displaystyle\int^{1}_{0}|\dot{q}_{R}(t)|^{2}dt}{\displaystyle\int^{1}_{0}(H-V(q_{R}(t)))dt},\label{21}
\end{eqnarray}
then by Lemmas 2.3-2.5,
$u_{R}(t)=q_{R}(\frac{t+\frac{T_{R}}{2}}{T_{R}}):\left(-\frac{T_{R}}{2},\frac{T_{R}}{2}\right)\rightarrow
\overline{\Lambda}_{R}$ is a non-constant $T_{R}$-periodic solution
satisfying (\ref{16}) and (\ref{17}).

\vspace{0.3cm}{\bf Remark 6}\ In our model, the set $\Lambda_{R}$ is
a closed set in the open set $E$. We minimize the functional on the
set $\Lambda_{R}$, however, we can not show that $u_{R}(t)$ solve
the equation at $\pm \frac{T_{R}}{2}$. But it is true that we do not
need that $u_{R}(t)$ is a solution at these two moments since
 we will take limits by letting $T_{R}\rightarrow+\infty$ later.
 Furthermore, we know that $u_{R}(t)$ still has definition at $\pm
 \frac{T_{R}}{2}$ and $|u_{R}(\pm \frac{T_{R}}{2})|=R$.

\vspace{0.3cm}{\bf Remark 7}\ The solution $u_{R}$ may have
collisions. If we need to prove that $u_{R}$ has no collision for
any $R>0$, as in the strong force case, one way is to prove the
potential satisfies the $Gordon's\ Strong\ Force$ condition. From
Remark 3, we can see $V$ satisfies $(B_{3})$ and $(B_{4})$ which are
the classical condition in strong force case, but in this paper, the
potential $V$ is positive, then we can not prove it satisfies the
$Gordon's\ Strong\ Force$ condition. However, in the following
lemmas, we can prove the minimizer has no collision which means it
is in $\Lambda_{R}$.

\section{Blowing-up Arguments}

\ \ \ \ \ \ In the following, we need to show that $u_{R}(t)$ can
not diverge to infinity uniformly as $R\rightarrow+\infty$.
Moreover, we prove the following lemma.

\vspace{0.3cm}{\bf Lemma 4.1}\ {\em Suppose that
$u_{R}(t):\left[-\frac{T_{R}}{2},\frac{T_{R}}{2}\right]\rightarrow
\overline{\Lambda}_{R}$ is the solution obtained in Lemma 3.1, then
$\min_{t\in\left[-\frac{T_{R}}{2},\frac{T_{R}}{2}\right]}|u_{R}(t)|$
is bounded from above . More precisely, there is a constant $M>0$
independent of $R$ such that
\begin{eqnarray*}
\min_{t\in\left[-\frac{T_{R}}{2},\frac{T_{R}}{2}\right]}|u_{R}(t)|\leq
M\ \ \ \mbox{for all}\ \ \ R>0.
\end{eqnarray*} }

\vspace{0.3cm}{\bf Proof.}\ Since $q_{R}\in M_{H}$, it is easy to
see that
\begin{eqnarray*}
\int^{\frac{T_{R}}{2}}_{-\frac{T_{R}}{2}}2H-(2V(u_{R}(t))+(\nabla
V(u_{R}(t)),u_{R}(t)))dt=0.
\end{eqnarray*}
There are two cases needed to be discussed.

{\bf Case 1.}\ $2H-(2V(u_{R}(t))+(\nabla
V(u_{R}(t)),u_{R}(t)))\equiv0$, which implies that
\begin{eqnarray*}
2H&=&2V(u_{R}(t))+(\nabla V(u_{R}(t)),u_{R}(t)),\ \ \ \mbox{a.e.} \
t\in\left[-\frac{T_{R}}{2},\frac{T_{R}}{2}\right].
\end{eqnarray*}
Remark 3 and hypothesis $(B_{3})$ imply that there exists a constant
$M_{1}>0$ independent of $R$ such that
\begin{eqnarray*}
\min_{t\in\left[-\frac{T_{R}}{2},\frac{T_{R}}{2}\right]}|u_{R}(t)|\leq
M_{1}.
\end{eqnarray*}

{\bf Case 2.}\ $2(H-V(u_{R}(t)))-(\nabla V(u_{R}(t)),u_{R}(t))$
changes sign in $\left[-\frac{T_{R}}{2},\frac{T_{R}}{2}\right]$.
Then there exists
$t_{0}\in\left[-\frac{T_{R}}{2},\frac{T_{R}}{2}\right]$ such that
\begin{eqnarray*}
2H-(2V(u_{R}(t_{0}))+(\nabla V(u_{R}(t_{0})),u_{R}(t_{0})))<0,
\end{eqnarray*}
which implies that
\begin{eqnarray*}
2H&<&2V(u_{R}(t_{0}))+(\nabla
V(u_{R}(t_{0})),u_{R}(t_{0}))=(2-\alpha)V(u_{R}(t_{0})).
\end{eqnarray*}
It follows from Remark 4 and hypothesis $(B_{6})$ that there exists
a constant $M_{2}>0$ independent of $R$ such that
\begin{eqnarray*}
\min_{t\in\left[-\frac{T_{R}}{2},\frac{T_{R}}{2}\right]}|u_{R}(t)|\leq
M_{2}.
\end{eqnarray*}
Then the proof is completed.

\section{Proof of Theorem 1.4}

\ \ \ \ \ \ The ideas for the following proofs in this section
mostly comes from Lemma 2.1 and Lemma 4.1 in \cite{3}, but there is
still some difference.

 \vspace{0.3cm}{\bf Lemma 5.1}\ {\em Suppose that $u_{R}(t)$
is the solution for $(1)-(2)$ on
$\left(-\frac{T_{R}}{2},\frac{T_{R}}{2}\right)$ obtained in Lemma
3.1. Then there exists a constant $m>0$ independent of $R$ such
that}
\begin{eqnarray*}
\min_{t\in\left[-\frac{T_{R}}{2},\frac{T_{R}}{2}\right]}|u_{R}(t)|\geq
m.
\end{eqnarray*}

\vspace{0.3cm}{\bf Proof.}\ Since $u_{R}(t)$ is a solution for
system $(1)-(2)$, then we can deduce that
\begin{eqnarray*}
\frac{d^{2}}{dt^{2}}\left(\frac{1}{2}|u_{R}(t)|^{2}\right)&=&\frac{d}{dt}(u_{R}(t),\dot{u}_{R}(t))\nonumber\\
&=& \mbox{} |\dot{u}_{R}(t)|^{2}+(u_{R}(t),\ddot{u}_{R}(t))\nonumber\\
&=& \mbox{} 2(H-V(u_{R}(t)))-(\nabla V(u_{R}(t)),u_{R}(t)),\ \ \ \
t\in\left(-\frac{T_{R}}{2},\frac{T_{R}}{2}\right).
\end{eqnarray*}
Since
$\left|u_{R}\left(-\frac{T_{R}}{2}\right)\right|=\left|u_{R}\left(\frac{T_{R}}{2}\right)\right|=R
$, then by Lemma 1.1 and hypothesis $(B_{4})$, we can find $m>0$
independent of $R$ such that
\begin{eqnarray*}
2H-(2V(u_{R}(t))+(\nabla V(u_{R}(t)),u_{R}(t)))<0,
\end{eqnarray*}
for any $t\in S=\{t\in[-\frac{T_{R}}{2},\frac{T_{R}}{2}]||u_{R}(t)|<
m\}$, which implies that $|u_{R}(t)|$ is concave when $t\in S$ and
we can deduce that $|u_{R}(t)|$ cannot take a local minimum in $S$,
which implies that $S=\emptyset$. If not, we can assume that there
exists a $\overline{t}\in S$, then we can easily check that
$|u_{R}(t)|$ takes a local minimum at some $\tilde{t}$ with
$|u_{R}(\tilde{t})|<m$, which is a contradiction. Then we have
\begin{eqnarray*}
 |u_{R}(t)|\geq m\ \ \ \mbox{for}\ \ \mbox{all}\ \ \
 t\in\left[-\frac{T_{R}}{2},\frac{T_{R}}{2}\right],
\end{eqnarray*}
which proves this lemma.

From this lemma, we can see that $u_{R}\in\Lambda_{R}$ has no
collision.

\vspace{0.3cm}{\bf Lemma 5.2}\ {\em Suppose that $R>M$,
$D\in(0,\min\{1,\frac{Hm^{\alpha}}{C_{2}}\})$ is a constant
independent of $R$ and $u_{R}(t)$ is the solution for
$(\ref{16})-(\ref{17})$ obtained in Lemma 3.1, where $M$ and $m$ are
from Lemma 4.1 and Lemma 5.1 respectively. Set
\begin{eqnarray*}
t_{+}=\sup\left\{t\in\left[-\frac{T_{R}}{2},\frac{T_{R}}{2}\right]|\left|u_{R}(t)\right|\leq
L\right\}
\end{eqnarray*}
and
\begin{eqnarray*}t_{-}=\inf\left\{t\in\left[-\frac{T_{R}}{2},\frac{T_{R}}{2}\right]|\left|u_{R}(t)\right|\leq
L\right\}
\end{eqnarray*}
where $L$ is a constant independent of $R$ such that $M<L<R$. Then
we have that}
\begin{eqnarray*}
\frac{T_{R}}{2}-t_{+}\rightarrow+\infty,\ \ \
t_{-}+\frac{T_{R}}{2}\rightarrow+\infty\ \ \ \mbox{as}\ \
R\rightarrow+\infty.
\end{eqnarray*}

\vspace{0.3cm}{\bf Proof.}\ By the definition of $u_{R}(t)$ we have
that
\begin{eqnarray*}
\left|u_{R}\left(-\frac{T_{R}}{2}\right)\right|=\left|u_{R}\left(\frac{T_{R}}{2}\right)\right|=R.
\end{eqnarray*}
Then, by Lemma 2.1 and the definitions of $t_{+}$ and $D$, we have
\begin{eqnarray}
\int^{\frac{T_{R}}{2}}_{t_{+}}\sqrt{H-DV(u_{R}(t))}|\dot{u}_{R}(t)|dt&\geq&\int^{\frac{T_{R}}{2}}_{t_{+}}\sqrt{H-\frac{DC_{2}}{|u_{R}(t)|^{\alpha}}}|\dot{u}_{R}(t)|dt\nonumber\\
&\geq& \mbox{}
\sqrt{H-\frac{DC_{2}}{m^{\alpha}}}\int^{\frac{T_{R}}{2}}_{t_{+}}|\dot{u}_{R}(t)|dt\nonumber\\
&\geq& \mbox{}
\sqrt{H-\frac{DC_{2}}{m^{\alpha}}}\left|\int^{t_{-}}_{-\frac{T_{R}}{2}}\dot{u}_{R}(t)dt\right|\nonumber\\
&\geq& \mbox{}\sqrt{H-\frac{DC_{2}}{m^{\alpha}}}(R-L).\label{3}
\end{eqnarray}
It follows from Lemma 2.1 and (\ref{17}) that
\begin{eqnarray*}
\int^{\frac{T_{R}}{2}}_{t_{+}}\sqrt{H-DV(u_{R}(t))}|\dot{u}_{R}(t)|dt&=&\sqrt{2}\int^{\frac{T_{R}}{2}}_{t_{+}}\sqrt{H-DV(u_{R}(t))}\sqrt{H-V(u_{R}(t))}dt\nonumber\\
&\leq& \mbox{} \sqrt{2}H\left(\frac{T_{R}}{2}-t_{+}\right)
\end{eqnarray*}
Combining (\ref{3}) with the  above estimation, we obtain that
\begin{eqnarray*}
\sqrt{H-\frac{DC_{2}}{m^{\alpha}}}(R-L)\leq\sqrt{2}H\left(\frac{T_{R}}{2}-t_{+}\right).
\end{eqnarray*}
Then we have
\begin{eqnarray*}
\frac{T_{R}}{2}-t_{+}\rightarrow+\infty,\ \ \ \mbox{as}\ \
R\rightarrow+\infty.
\end{eqnarray*}
The limit for $t_{-}+\frac{T_{R}}{2}$ can be obtained in the similar
way. The proof is completed.

Subsequently, we set that
\begin{eqnarray}
t^{*}=\inf\left\{t\in\left[-\frac{T_{R}}{2},\frac{T_{R}}{2}\right]||u_{R}(t)|=M\right\}\label{31}
\end{eqnarray}
and
\begin{eqnarray*}
u_{R}^{*}(t)=u_{R}(t-t^{*})
\end{eqnarray*}
Since $L>M$, we can deduce that $t_{+}\geq t^{*}\geq t_{-}$, which
implies that
\begin{eqnarray*}
-\frac{T_{R}}{2}+t^{*}\rightarrow-\infty,\ \ \mbox{}\ \
\frac{T_{R}}{2}+t^{*}\rightarrow+\infty\ \ \mbox{as}\ \
R\rightarrow\infty.
\end{eqnarray*}

It follows from (\ref{17}) that
\begin{eqnarray*}
   \frac{1}{2}|\dot{u}^{*}_{R}(t)|^{2}+ V(u^{*}_{R}(t))=H,\ \ \
   \forall\ t\in\left(-\frac{T_{R}}{2}+t^{*},\frac{T_{R}}{2}+t^{*}\right).
\end{eqnarray*}
which implies that
\begin{eqnarray*}
|\dot{u}^{*}_{R}(t)|^{2}=2(H-V(u^{*}_{R}(t))),\ \ \
  \forall\ t\in\left(-\frac{T_{R}}{2}+t^{*},\frac{T_{R}}{2}+t^{*}\right).
\end{eqnarray*}
By Lemma 5.1,  Remark 3, $(B_{3})$ and $V\in
C^{1}(R^{N}\setminus\{0\},R^{1})$, we can deduce that there exists a
constant $M_{4}>0$ independent of $R$ such that
\begin{eqnarray*}
 |V(u^{*}_{R}(t))|\leq M_{4}\ \ \ \mbox{for all}\ \
 t\in\left[-\frac{T_{R}}{2}+t^{*},\frac{T_{R}}{2}+t^{*}\right].
\end{eqnarray*}
Then there is a constant $M_{5}$ independent of $R$ such that
\begin{eqnarray*}
 |\dot{u}^{*}_{R}(t)|\leq M_{5}\ \ \ \mbox{for all}\ \
 t\in\left[-\frac{T_{R}}{2}+t^{*},\frac{T_{R}}{2}+t^{*}\right].
\end{eqnarray*}
which implies that
\begin{eqnarray*}
\ |u_{R}^{*}(t_{1})-u_{R}^{*}(t_{2})| \leq
\left|\int^{t_{1}}_{t_{2}}\dot{u}_{R}^{*}(s)ds\right| \leq \mbox{}
\int^{t_{1}}_{t_{2}}\left|\dot{u}_{R}^{*}(s)\right|ds \leq
M_{5}|t_{1}-t_{2}|
\end{eqnarray*}
for each $R>0$ and $t_{1}, t_{2} \in
\left[-\frac{T_{R}}{2}+t^{*},\frac{T_{R}}{2}+t^{*}\right]$, which
shows $\{u_{R}^{*}\}$ is equicontinuous. Then there is a
subsequence $\{u_{R}^{*}\}_{R>0}$ converging to $u_{\infty}$ in
$C_{loc}(R^{1},R^{N})$. Then there exists a function
$u_{\infty}(t)$ such that
\begin{eqnarray*}
&&(\mbox{i})\ u_{R}^{*}(t)\rightarrow u_{\infty}(t)\ \mbox{in}\ C_{loc}(R^{1},R^{N})\nonumber\\
&&(\mbox{ii})|u_{\infty}(t)|\rightarrow +\infty\ \mbox{as}\
|t|\rightarrow+\infty
\end{eqnarray*}
and $u_{\infty}(t)$ satisfies systems $(1)-(2)$.

From the above lemmas, we have proved there is at least one
hyperbolic solution for $(1)-(2)$ with $H>0$.

\section{Proof of Theorem 1.5}
\ \ \ \ \ \ By the conditions of Theorem 1.5, the existence of
hyperbolic solutions for systems $(1)-(2)$ can be obtained with a
similar proof of Theorem 1.4. Subsequently, we give the proof of the
asymptotic direction of hyperbolic solutions at infinity. The proof
is similar to Felmer and Tanaka's in \cite{3}.

\vspace{0.3cm}{\bf Lemma 6.1}\ {\em Suppose that $u_{R}(t)$ is the
solution for $(\ref{16})-(\ref{17})$ obtained in Lemma 3.1. Then
there exists a constant $M_{6}>0$ independent of $R>1$ such that}
\begin{eqnarray*}
\int^{\frac{T_{R}}{2}}_{-\frac{T_{R}}{2}}\sqrt{H-V(u_{R}(t))}|\dot{u}_{R}(t)|dt\leq
\sqrt{2H}R+M_{6}.
\end{eqnarray*}

\vspace{0.3cm}{\bf Proof.}\ Firstly, we define the function $\xi(t)$
on $[1,+\infty)$ as a solution of
\begin{eqnarray*}
&& \dot{\xi}(t)=\sqrt{2(H-V(\xi(t)e))}\\
&& \xi(1)=1.
\end{eqnarray*}
And $\tau_{R}>1$ is a real number such that $\xi(\tau_{R})=R$. We
can define $\xi(t)$ in $(-\infty,0]$ and $\tau_{-R}$ in a similar
way. Then we can fix $\varphi(t)\in H^{1}([0,1],R^{N})$ such that
$\tilde{\gamma}_{R}(t)\in \Lambda_{R}$ where
\begin{eqnarray*}
\tilde{\gamma}_{R}(t)=\gamma_{R}(t(\tau_{R}-\tau_{-R})+\tau_{-R}),\
\ \ \mbox{and} \ \ \gamma_{R}(t)=\left\{
\begin{array}{ll}
\xi(t)e&\mbox{for $t\in[1,\tau_{R}]\bigcup[\tau_{-R},0]$}\\
\varphi(t)&\mbox{for $t\in[0,1]$}.
\end{array}
\right.
\end{eqnarray*}
Subsequently, we set $u_{r}(t)=\tilde{\gamma}_{R}(\frac{t+r}{2r})$.
And it is easy to see that $u_{r}(t)=\gamma_{R}(t)$ if $\pm
r=\tau_{\pm R}$. Similar to \cite{3}, we can deduce that for $r>0$
\begin{eqnarray}
(2f(\tilde{\gamma}_{R}))^{\frac{1}{2}}&=&\inf_{r>0}\frac{1}{\sqrt{2}}\int^{r}_{-r}\frac{1}{2}|\dot{u}_{r}(t)|^{2}+H-V(u_{r}(t))dt\nonumber\\
&\leq& \mbox{}
\frac{1}{\sqrt{2}}\int^{\tau_{R}}_{-\tau_{R}}\frac{1}{2}|\dot{\gamma}_{R}(t)|^{2}+H-V(\gamma_{R}(t))dt.\label{26}
\end{eqnarray}
Since
$[-\tau_{R},\tau_{R}]=[-\tau_{R},0]\bigcup[0,1]\bigcup[1,\tau_{R}]$,
then by $(V_{1})$, we can estimate (\ref{26}) by three integrals.
Firstly, we estimate the integral on $[1,\tau_{R}]$
\begin{eqnarray*}
I_{[1,\tau_{R}]}&=&\frac{1}{\sqrt{2}}\int^{\tau_{R}}_{1}\frac{1}{2}|\dot{\gamma}_{R}(t)|^{2}+H-V(\gamma_{R}(t))dt\nonumber\\
&=& \mbox{}\int^{\tau_{R}}_{1}\sqrt{H-V(\xi(t)e)}\dot{\xi}(t)dt\nonumber\\
&=& \mbox{} \int^{R}_{1} \sqrt{H-V(se)}ds\nonumber\\
&\leq& \mbox{}  \sqrt{H}R+M_{7}
\end{eqnarray*}
for some $M_{7}>0$ independent of $R>1$. Similarly, we can get
\begin{eqnarray*}
I_{[\tau_{R},0]}\leq \sqrt{H}R + M_{7}.
\end{eqnarray*}
Since $I_{[0,1]}$ is independent of $R$, we obtain that
\begin{eqnarray*}
\frac{1}{\sqrt{2}}\int^{\tau_{R}}_{-\tau_{R}}\frac{1}{2}|\dot{\gamma}_{R}(t)|^{2}+H-V(\gamma_{R}(t))dt\leq
2\sqrt{H}R + M_{6}
\end{eqnarray*}
for some $M_{6}>0$ independent of $R$. Since $q(t)$ is the minimizer
of $f$ on $\Lambda_{R}$, then by (\ref{26}), we have
\begin{eqnarray}
\int^{\frac{T_{R}}{2}}_{-\frac{T_{R}}{2}}\sqrt{H-V(u_{R}(t))}|\dot{u}_{R}(t)|dt
&\leq&
\left(\int^{\frac{T_{R}}{2}}_{-\frac{T_{R}}{2}}H-V(u_{R}(t))dt\right)^{\frac{1}{2}}\left(\int^{\frac{T_{R}}{2}}_{-\frac{T_{R}}{2}}|\dot{u}_{R}(t)|^{2}dt\right)^{\frac{1}{2}}
\nonumber\\
&=&\mbox{}(2f(q))^{\frac{1}{2}}\nonumber\\
&\leq&\mbox{}(2f(\tilde{\gamma}_{R}))^{\frac{1}{2}}\nonumber\\
&\leq& \mbox{}
\frac{1}{\sqrt{2}}\int^{\tau_{R}}_{-\tau_{R}}\frac{1}{2}|\dot{\gamma}_{R}(t)|^{2}+H-V(\gamma_{R}(t))dt\nonumber\\
&\leq& \mbox{}2\sqrt{H}R + M_{6}.
\end{eqnarray}
Then we finish the proof of this lemma.

Similar to Felmer and Tanaka \cite{3}, we set
\begin{eqnarray*}
A(t)=\sqrt{|u_{R}(t)|^{2}|\dot{u}_{R}(t)|^{2}-(u_{R}(t),\dot{u}_{R}(t))^{2}}
\end{eqnarray*}
and
\begin{eqnarray*}
\omega(t)=\frac{A(t)}{|u_{R}(t)||\dot{u}_{R}(t)|}.
\end{eqnarray*}
Using the motion and energy equations, we have
\begin{eqnarray*}
|\dot{A}(t)|\leq |u_{R}(t)||\nabla V(u_{R}(t))|
\end{eqnarray*}
and
\begin{eqnarray*}
\frac{d\omega}{dt}=\frac{2}{|u_{R}(t)||\dot{u}_{R}(t)|}(-\omega\sqrt{1-\omega^{2}}sign(u_{R}(t),\dot{u}_{R}(t))(H-V(u_{R}(t)))+|u_{R}(t)||\nabla
V(u_{R}(t))|).
\end{eqnarray*}
The proof of the following lemma is the same as \cite{3}.

\vspace{0.3cm}{\bf Lemma 6.2(See\cite{3})}\ {\em Assume $u_{R}$ is a
solution for $(\ref{16})-(\ref{17})$ obtained in Lemma 3.1. For any
$\eta\in(0,1)$, there exists a $L_{\eta}\geq m$ such that if}
\begin{eqnarray}
|u_{R}(t_{0})|\geq L_{\eta},\ \ \
(u_{R}(t_{0}),\dot{u}_{R}(t_{0}))>0\ \ \ \mbox{and}\ \ \
\omega(t_{0})<\eta\label{28}
\end{eqnarray}
{\em for some $t_{0}\in (-\frac{T_{R}}{2},\frac{T_{R}}{2})$, then we
have for $t\in[t_{0},\frac{T_{R}}{2}]$}
\begin{eqnarray*}
&&(\mbox{i}).\ \omega(t)<\eta,\nonumber\\
&&(\mbox{ii}).\ \frac{d}{dt}|u_{R}(t)|\geq\sqrt{1-\eta^{2}}|\dot{u}_{R}(t)|,\nonumber\\
&&(\mbox{iii}).\ \frac{d}{dt}|u_{R}(t)|\geq\sqrt{2(1-\eta^{2})H},\nonumber\\
&&(\mbox{iv}).\
|u_{R}(t)|\geq|u_{R}(t_{0})|+\sqrt{2(1-\eta^{2})H}(t-t_{0}).
\end{eqnarray*}

\vspace{0.3cm}{\bf Lemma 6.3(See\cite{3})}\ {\em Let $u_{R}$ is a
solution for $(\ref{16})-(\ref{17})$ obtained in Lemma 3.1
satisfying (\ref{28}) and $|u_{R}(t)|\geq r_{0}$ with $t\geq t_{0}$
for certain $t_{0}\in (-\frac{T_{R}}{2},\frac{T_{R}}{2})$ with
$\eta\in(0,\frac{1}{2})$ and $L_{\eta}$ as in Lemma 6.2. Then for
$t\geq t_{0}$ we have}
\begin{eqnarray*}
\left|\frac{u_{R}(t)}{|u_{R}(t)|}-\frac{u_{R}(t_{0})}{|u_{R}(t_{0})|}\right|\leq
M_{8}\eta+\frac{M_{9}}{|u_{R}(t_{0})|^{\beta}},
\end{eqnarray*}
where $M_{8}$, $M_{9}>0$ are independent of $\eta$, $u_{R}(t)$ and
$t_{0}$.

\vspace{0.3cm}{\bf Proof.}\ By Lemma 5.1, (iii) of Lemma 6.2 and
$(V_{2})$, we can estimate $A(t)$ as following.
\begin{eqnarray}
A(t)&\leq& A(t_{0})+\int^{t}_{t_{0}}|u_{R}(s)||\nabla V(u_{R}(s))|ds\nonumber\\
&\leq& \mbox{}
A(t_{0})+\frac{M_{9}}{\sqrt{2(1-\eta^{2})H}}\int^{t}_{t_{0}}|u_{R}(s)||\nabla
V(u_{R}(s))|\frac{d}{ds}|u_{R}(s)|ds\nonumber\\
&\leq& \mbox{}
A(t_{0})+\frac{M_{9}}{\sqrt{2(1-\eta^{2})H}}\int^{|u_{R}(t)|}_{|u_{R}(t_{0})|}\varphi\left|\nabla
V\left( \frac{u_{R}(s)}{|u_{R}(s)|}\varphi\right)\right|d\varphi\nonumber\\
&\leq& \mbox{}
A(t_{0})+\frac{M_{9}M_{0}}{\sqrt{2(1-\eta^{2})H}}\int^{|u_{R}(t)|}_{|u_{R}(t_{0})|}\frac{1}{\varphi^{\beta}}d\varphi\nonumber\\
&\leq& \mbox{}
A(t_{0})+\frac{M_{9}M_{0}}{\sqrt{2(1-\eta^{2})H}(\beta-1)}\frac{1}{|u_{R}(t_{0})|^{\beta-1}}\nonumber\\
&\leq& \mbox{}A(t_{0})+
\frac{M_{10}}{|u_{R}(t_{0})|^{\beta-1}}\label{29}
\end{eqnarray}
for some $M_{10}>0$ independent of $R$. Since we have
\begin{eqnarray}
\left|\frac{d}{dt}\frac{u_{R}(t)}{|u_{R}(t)|}\right|=\frac{A(t)}{|u_{R}(t)|^{2}},\label{30}
\end{eqnarray}
then it follows from (iii) of Lemma 6.2, (\ref{29}) and (\ref{30})
that
\begin{eqnarray*}
\left|\frac{u_{R}(t)}{|u_{R}(t)|}-\frac{u_{R}(t_{0})}{|u_{R}(t_{0})|}\right|&\leq&\int^{t}_{t_{0}}\frac{A(s)}{|u_{R}(s)|^{2}}ds\nonumber\\
&\leq& \mbox{}\left(A(t_{0})+
\frac{M_{10}}{|u_{R}(t_{0})|^{\beta-1}}\right)\int^{t}_{t_{0}}\frac{1}{|u_{R}(s)|^{2}}ds\nonumber\\
&\leq& \mbox{}\left(A(t_{0})+
\frac{M_{10}}{|u_{R}(t_{0})|^{\beta-1}}\right)\frac{1}{\sqrt{2(1-\eta^{2})H}}\int^{t}_{t_{0}}\frac{1}{|u_{R}(s)|^{2}}\frac{d}{ds}|u_{R}(s)|ds\nonumber\\
&\leq& \mbox{}\left(A(t_{0})+
\frac{M_{10}}{|u_{R}(t_{0})|^{\beta-1}}\right)\frac{1}{\sqrt{2(1-\eta^{2})H}}\frac{1}{|u_{R}(t_{0})|}.
\end{eqnarray*}
By energy equation and the definition of $t_{0}$, we have
\begin{eqnarray*}
A(t_{0})=\omega(t_{0})|u_{R}(t_{0})||\dot{u}_{R}(t_{0})|\leq
\eta|u_{R}(t_{0})|\sqrt{2(H-V(u_{R}(t_{0})))},
\end{eqnarray*}
which implies that for some $M_{8}$, $M_{9}>0$ independent of $R$
\begin{eqnarray*}
\left|\frac{u_{R}(t)}{|u_{R}(t)|}-\frac{u_{R}(t_{0})}{|u_{R}(t_{0})|}\right|\leq
M_{8}\eta+\frac{M_{9}}{|u_{R}(t_{0})|^{\beta}},
\end{eqnarray*}
which proves this lemma.

Since we have Theorems 6.1-6.3, similar to \cite{3}, we have the
following theorem.

\vspace{0.3cm}{\bf Lemma 6.4(See\cite{3})}\ {\em For any
$\varepsilon>0$, there exists $M_{11}>0$ such that for $R>M_{11}$}
\begin{eqnarray*}
u_{R}\left(\left[t^{*},\frac{T_{R}}{2}\right]\right)\bigcap\{|x|\geq
M_{11}\}\subset\left\{y\in
R^{N}:\left|\frac{y}{|y|}-e\right|<\varepsilon\right\},
\end{eqnarray*}
where $e$ is the given direction defined in $\Lambda_{R}$ and
$t^{*}$ is defined as (\ref{31}).

Let $\overline{t}\geq t^{*}$ such that
$|u_{R}(\overline{t})|=\overline{L}_{\eta}$. Then we can get for any
$\varepsilon>0$
\begin{eqnarray}
\left|\frac{u_{R}(t)}{|u_{R}(t)|}-e\right|<\varepsilon
\end{eqnarray}
for all $t\geq\bar{t}$, which implies that
\begin{eqnarray*}
\frac{u_{\infty}(t)}{|u_{\infty}(t)|}\rightarrow e\ \ \ \mbox{as}\ \
\ t\rightarrow+\infty
\end{eqnarray*}
and
\begin{eqnarray*}
\frac{u_{\infty}(t)}{|u_{\infty}(t)|}\rightarrow e\ \ \ \mbox{as}\ \
\ t\rightarrow-\infty.
\end{eqnarray*}

From the above discussion, we have proved there is at least one
hyperbolic solution for $(1)-(2)$ with $H>0$ which has the given
asymptotic direction at infinity. We finish the proof. $\Box$

\end{document}